\documentclass[12pt]{article}
\usepackage{amsmath}
\usepackage{amssymb}

\newcommand{\llim}{\lim\limits}
\newcommand{\lint}{\int\limits}
\newcommand{\qed}{\mbox{$\quad\blacksquare$}}
\begin{document}
\hspace{-2cm}
\raisebox{12ex}[1ex]{\fbox{{\footnotesize A 
slightly different version will appear 
in {\it American Mathematical Monthly}}}}
\begin{center}
{\large\bf Some divergent trigonometric integrals}
\vskip.25in
Erik Talvila\\ [2mm]
{\footnotesize
Department of Mathematical Sciences\\
University of Alberta\\
Edmonton AB Canada T6G 2E2\\
etalvila@math.ualberta.ca\\
November 22, 2000}
\end{center}
{\bf 1. Introduction.}\quad
Browsing through an integral table on a dull Sunday afternoon some
time ago, I came across four divergent trigonometric integrals.
(See \eqref{I.1} and \eqref{I.2} below.)  I was intrigued as to
how these divergent integrals ended up in a respectable table.  Tracing
their history, it turned out they were originally ``evaluated''
when some convergent integrals, \eqref{I.4} and \eqref{I.5},
were differentiated under the integral sign with respect to a parameter,
formally yielding \eqref{I.1} and \eqref{I.2}.  We will give a simple
proof that these integrals diverge, look at their history in print and then
make some final remarks about necessary and sufficient conditions for 
differentiating
under the integral sign.  We have no
motive
in defaming either 
the (shockingly famous) mathematician who made the original error,
or the editors of the otherwise fine tables in which the integrals appear.
We all make mistakes and we're not out to  point the finger at anyone.  (In 
this regard see the last two exercises of Chapter~2 in
\cite{spivak}.)  
We will also see that Maple and Mathematica have considerable difficulties
with these integrals.

\bigskip
\noindent
{\bf 2. Four divergent integrals.}\quad
Here they are.  Throughout, $a$ and $b$ are positive real numbers.  Purported
values appear on the right:
\begin{eqnarray}
\lint_{x=0}^{\infty}x\left\{\!\!\begin{array}{c}
\sin(ax^2)\\
\cos(ax^2)
\end{array}
\!\!\right\}
\sin(bx)\,dx & \mbox{`}=\mbox{'} & 
\frac{b}{4a}\sqrt{\frac{\pi}{2a}}\left[\sin\left(\frac{b^2}{4a}\right)
\pm\cos\left(\frac{b^2}{4a}\right)\right]\label{I.1}\\
\lint_{x=0}^{\infty}x\left\{\!\!\begin{array}{c}
\sin(ax^2)\\ 
\cos(ax^2)
\end{array}
\!\!\right\}
\cos(bx)\,dx & \mbox{`}=\mbox{'} & \label{I.2} 
\end{eqnarray}
\vspace{-5\jot}
$$
\frac{1}{2a}\left\{\!\!\begin{array}{c}
1\\
0
\end{array}\!\!\right\}\mp\frac{b}{2a}\sqrt{\frac{\pi}{2a}}
\left[\left\{\!\!\begin{array}{c}\sin\left[b^2/(4a)\right]\\
\cos\left[b^2/(4a)\right]
\end{array}
\!\!\right\}
\rm{C}\left(\frac{b^2}{4a}\right)\mp
\left\{\!\!\begin{array}{c}\cos\left[b^2/(4a)\right]\\
\sin\left[b^2/(4a)\right]
\!\!\end{array}
\right\}
\rm{S}\left(\frac{b^2}{4a}\right)\right].
$$
The two Fresnel integrals are
\begin{equation}
C(x)= 
\frac{1}{\sqrt{2\pi}}\int_{0}^x\!\!\cos t\frac{dt}{\sqrt{t}}\quad\mbox{ and }
\quad 
S(x)=\frac{1}{\sqrt{2\pi}}\int_{0}^x\!\!\sin t\frac{dt}{\sqrt{t}}\,.
\end{equation}
Note that in the literature, the same symbols $C$ and $S$ denote 
several different 
definitions of the Fresnel integrals.

Let's prove  these devils diverge.  \\

\noindent
{\bf Proposition:}\quad
{\em The integrals in \eqref{I.1} and \eqref{I.2} diverge}.

\bigskip
{\bf Proof:} Consider $A:= \int_{-\infty}^\infty x e^{i(x^2+x)}\,dx$ 
(which, unfortunately, does not exist).  Since the integrand is continuous,
this integral exists if and only if the limits
$$
\llim_{T\to\infty}\int_0^T x e^{i(x^2+x)}\,dx\qquad
\llim_{T\to\infty}\int^0_{-T} x e^{i(x^2+x)}\,dx 
$$
both exist.
Let $T_1$, $T_2>0$.  Integrate 
by parts and complete the square:
\begin{equation}
\int^{T_2}_{-T_1} x e^{i(x^2+x)}\,dx = \frac{1}{2i}\left[e^{i(T_2^2+T_2)}
-e^{i(T_1^2-T_1)}\right]-\frac{e^{-i/4}}{2}\lint_{-T_1+1/2}^{T_2+1/2}e^{ix^2}\,dx.
\label{I.3}
\end{equation}
Let's look at the convergence of $I:=\int_{-\infty}^\infty e^{ix^2}\,dx$.
Use the substitution $x^2=t$.  Then $I=\int_0^\infty e^{it}dt/\sqrt{t}=
\sqrt{2\pi}(C(\infty)+iS(\infty))$.  This
can be seen to converge by applying Dirichlet's Test  over $(1,\infty)$. 
(For Dirichlet's Test see, for example, \cite{rogers}, pp~261.)
More properly, we should start with the `$t$' version of $I$, show it converges
and then transform back to the `$x$' version.  But let's keep it cool.  In fact,
many ways have been found to evaluate $I$.   One method is to use contour 
integration. (Rotate the integral $\int_{-\infty}^\infty e^{-x^2}dx$ by $\pi/4$ in
the complex plane.  See \cite{spiegel}, pp~184.) A second method is to use
the gamma function \cite{rogers}, pp~272.
The result is $I=e^{i\pi/4}\sqrt{\pi}$.  Now, as $T_1$, $T_2\to\infty$ (independently
of each other) the final integral in \eqref{I.3} becomes $I$ but the bracketed
term fails to have a limit.  Hence, the integral $A$ diverges.

To get the integrals in \eqref{I.1} and \eqref{I.2} we do the following.
Suppose $B:=\int_{-\infty}^\infty x e^{i(x^2-x)}\,dx$ converged.  The
transformation $x\mapsto -x$ gives $B=-A$ so $B$ diverges.  The transformations
$x\mapsto x/\sqrt{a}\pm (\sqrt{a}-b)/(2a)$ now show that
$C:=\int_{-\infty}^\infty x e^{i(ax^2\pm bx)}\,dx$ diverges for all positive $a$ 
and $b$.  Finally, if the integrals in \eqref{I.1} and \eqref{I.2}
converged then we could form the four convergent linear combinations
$$
\lint_{x=0}^{\infty}x\left[\cos(ax^2)\cos(bx)\mp\sin(ax^2)\sin(bx)\right]
\,dx
$$
$$
\lint_{x=0}^{\infty}x\left[\sin(ax^2)\cos(bx)\pm\cos(ax^2)\sin(bx)\right]
\,dx.
$$
But, using the addition formulas for the sine
and cosine functions and then Euler's formula yields $C$.  Hence, 
the integrals in \eqref{I.1} and \eqref{I.2} diverge.
\qed

To see the manner in which the integrals diverge, let
\begin{eqnarray*}
A_T & := & \int_{-T}^T x e^{i(x^2+x)}\,dx\\
 & = & e^{iT^2}\sin T-\frac{e^{-i/4}}{2}\int_{-T+1/2}^{T+1/2}  e^{ix^2}\,dx.
\end{eqnarray*}
As $T\to\infty$, the integral term in the line 
above has limit $\sqrt{\pi}\,e^{i(\pi-1)/4}
/2$, whereas the term $e^{iT^2}\sin T$ oscillates rapidly with unit
magnitude.
Note that this
also shows that our integrals do not even exist as Cauchy principal
values.

\bigskip
\noindent
{\bf 3. History of the divergent integrals.}\quad
Now we'll look at the history of \eqref{I.1} and \eqref{I.2} in print.
The thickest book of integrals (3500 pages in five volumes) is that
of Prudnikov\footnote[2]{Recently deceased.  A heartfelt obituary by
Marichev appears in \cite{marichev}.  See also \cite{glaeske}.},
Brychkov and Marichev, \cite{prudnikov2}. Our
integrals appear in Volume I,  2.5.22.  They also appear in the original Russian edition
\cite{prudnikov1}.  Sources are not
referenced in this work.  It is interesting that they do not appear in the earlier
book \cite{ditkin} by Ditkin and Prudnikov.  Other major tables they are absent
from are \cite{erdelyi},
\cite{grobner}, \cite{ober1} and \cite{ober2}.  As the tables by Erd\'{e}lyi
and Oberhettinger are quite comprehensive, one suspects it was noticed that these
integrals diverged and they were purposely omitted.  However, they are contained in the 
Gradshteyn and Ryzhik tome,
\cite{gradshteyn} (3.851).
They are not in the first few Russian editions but the 1963 edition \cite{grad1}
was enlarged enough to include these divergent integrals.
All subsequent Russian editions and 
English
translations beginning 1965 \cite{grad2} 
contain the integrals in \eqref{I.1} and \eqref{I.2}.
Now, Gradshteyn and Ryzhik do give references.  They say (in a garbled
citation) that our
integrals come from tables by Bierens de Haan.  

David Bierens de Haan (1822--1895) was a Dutch mathematician noted for
compiling tables of integrals, for actuarial work, for writing various
essays in the history of science and mathematics, for producing an
encyclop{\ae}dic biography of Dutch scientists and for being an early editor
of the works of Christian Huygens. (A mammoth task.  It took until 1950 when
the 22$^{nd}$ volume was finally published.)
A complete list of de Haan's publications is given in
\cite{korteweg}.  There have been several papers on his
life and work.    See \cite{schrek} for references,
photos and an interesting
reproduction of the 1935 title page from a Japanese edition of his integral
table.  His 1858 {\it Tables d'int\'egrales d\'efinies} \cite{haan1858} was 
the first really substantial
table of integrals.
It was enlarged and corrected in
an 1867 edition \cite{haan1867}.  For nearly a century these were the preeminent
integral tables.  The 1867 edition was still being  reprinted in  1957 
\cite{haan1957}, three years after the publication of the Bateman Manuscript
tables \cite{erdelyi}.  The integrals in \eqref{I.1} appear in the 1858 table
\cite{haan1858}, formulas 193.17 and 193.18, an 1862 companion 
volume that details the techniques used
to compute integrals in the tables \cite{haan1862}, pp~443, and in the 1867
table \cite{haan1867}, formulas 150.4, 150.7.  

Now, Bierens de Haan lists Cauchy as his source for \eqref{I.1}.  An
examination of Cauchy's works (27 volumes!) shows these integrals appear
twice \cite{cauchy1825} (1815) and \cite{cauchy1827} (see also \cite{cauchyger})
(1825).  In both instances, 
Cauchy correctly obtains the convergent integrals
\begin{equation}
\lint_{x=0}^{\infty}\left\{\!\!\begin{array}{c}
\sin(ax^2)\\
\cos(ax^2)
\end{array}
\!\!\right\}
\cos(bx)\,dx  = 
\frac{1}{2}\sqrt{\frac{\pi}{2a}}\left[\cos\left(\frac{b^2}{4a}\right)
\mp\sin\left(\frac{b^2}{4a}\right)\right]\label{I.4}
\end{equation}
with $a=1$.  As is clearly stated in the above references, he then proceeds
to differentiate under the integral sign with respect to $b$.  Very bad!
The functions defined by the integrals in \eqref{I.4} are certainly
differentiable since the right side of \eqref{I.4} is differentiable.  But,
differentiating under the integral sign leads to our divergent integrals
\eqref{I.1}.  It is not some unknown schmo but Cauchy, the ``Father of rigour", 
who commits an 
error that has been copied for 185 years.  (The appropriateness of
this epithet is contested.  One in favour is \cite{grabiner}.
One against is \cite{gratton}.)  Bierens de Haan repeats
this argument in \cite{haan1862}, pp~443.

When the two integrals
\begin{eqnarray}
\lefteqn{\lint_{x=0}^{\infty}\left\{\!\!\begin{array}{c}
\sin(ax^2)\\
\cos(ax^2)
\end{array}
\!\!\right\}
\sin(bx)\,dx  =} \notag\\
 & & \sqrt{\frac{\pi}{2a}}
\left[\left\{\!\!\begin{array}{c}\cos\left[b^2/(4a)\right]\\
\sin\left[b^2/(4a)\right]
\end{array}
\!\!\right\}
\rm{C}\left(\frac{b^2}{4a}\right)\pm
\left\{\!\!\begin{array}{c}\sin\left[b^2/(4a)\right]\\
\cos\left[b^2/(4a)\right]
\!\!\end{array}
\right\}
\rm{S}\left(\frac{b^2}{4a}\right)\right]
\label{I.5}
\end{eqnarray}
are differentiated under the integral sign we get the divergent integrals
\eqref{I.2}.  After some incorrect manipulations, Bierens de Haan obtains the value
$0$ for these integrals \cite{haan1862}, pp~443.  He then differentiates under the
integral sign to get $0$ for \eqref{I.2}.  

The integrals in \eqref{I.4} and \eqref{I.5} may be evaluated using the methods
in the proof of the Proposition.

The tables of Bierens de Haan had many errors, both mathematical and 
typographical.  Two long works discussing the correctness of
his tables are \cite{lindman} and
\cite{sheldon}.  Neither mentions our divergent integrals.
All of the integral tables listed above have received considerable scrutiny.
The journal, {\it Mathematics of computation}, and its predecessor,
{\it Mathematical tables and other aids to computation}, list numerous errata.
However, despite over 300 published pages of errata related to the above
tables there do not seem to be
any references to 
\eqref{I.1} and \eqref{I.2}.  The article \cite{klerer} compares the
correctness of
various integral tables.
The interested reader should consult this article
to see how shockingly high the error rates are.

\bigskip
\noindent
{\bf 4. Maple and Mathematica.}\quad
Here are how Maple (V.4) and Mathematica (4.0) fare.  Maple correctly evaluates 
\eqref{I.4} and \eqref{I.5} for arbitrary $a$ and $b$ but falters when asked to
perform the calculation with specific numerical values.
For example, it gives $\int_0^\infty \sin(3.1x^2)\cos(2.2x)\,dx=0$.  Maple
correctly says the integrals \eqref{I.1} and \eqref{I.2} diverge.  Mathematica
fails in a different way.  It correctly calculates \eqref{I.4} and \eqref{I.5}
(considerable simplification needed to obtain the form of \eqref{I.5} above).
But, Mathematica thinks \eqref{I.1} and \eqref{I.2} converge!  It gives the
same incorrect values that are in the tables.

\bigskip
\noindent
{\bf 5. Differentiation under the integral sign.}\quad
Differentiating the convergent integrals \eqref{I.4} and \eqref{I.5} under
the integral sign with respect to $b$ yields the divergent integrals
\eqref{I.1} and \eqref{I.2}.  This doesn't mean the functions defined by
\eqref{I.4} and \eqref{I.5} aren't differentiable it just means we cannot
obtain their derivatives by differentiating under the integral.  Could we
have predicted this in advance?  This is a difficult problem.  Suppose
we have $\int_a^bf(x,y)\,dy$.  A sufficient condition for differentiating under
the integral, suitable for Riemann integrals,
is that $\int_a^bf_1(x,y)\,dy$ converge uniformly in $x$.
See \cite{rogers}, pp~260.  For Lebesgue
integrals the dominating condition $|f_1(x,y)|\leq g(y)$ where $g\in L^1$
suffices.
For Riemann and Lebesgue integrals, necessary and sufficient conditions for
differentiating under the
integral sign  are harder to come by.  However, this is a much simpler
problem when we use the Henstock integral.  The solution depends on
being able to integrate every derivative, a property not held by
either the Riemann or the Lebesgue integral.  The interested reader
can see \cite{talvila}.
A good introduction to
the Henstock integral is given in \cite{bartle}.

\end{document}